\documentclass[runningheads]{llncs}
\usepackage[T1]{fontenc}
\usepackage{graphicx}
\usepackage{amsmath}
\usepackage{amsfonts}
\usepackage{gensymb}
\usepackage{url}
\usepackage{pseudo}

\newcommand{\bsx}{\boldsymbol x}

\newcommand{\bsi}{\boldsymbol I}

\newcommand{\bsD}{\boldsymbol\Delta}
\newcommand{\bsbe}{\boldsymbol b}

\newcommand{\IR}{\mathbb R}
\newcommand{\calR}{\mathcal R}

\begin{document}

\title{Pushing the Limits of the Reactive Affine Shaker Algorithm to Higher Dimensions}
\titlerunning{Pushing the Limits of the Reactive Affine Shaker to Higher Dimensions}
\author{Roberto Battiti\orcidID{0000-0002-0259-8603} \and
Mauro Brunato\orcidID{0000-0002-7885-4255}}
\institute{DISI --- Università di Trento, Via Sommarive 9, I-38123 Trento, Italy\\
\email{\{roberto.battiti,mauro.brunato\}@unitn.it}}

\maketitle

\begin{abstract}

Bayesian Optimization (BO) for the minimization of expensive functions of continuous variables 
uses all the knowledge acquired from previous samples ($\bsx_i$ and $f(\bsx_i)$ values) to build
a surrogate model based on Gaussian processes. The surrogate is then exploited to define the next point to sample,
through a careful balance of exploration and exploitation.
Initially intended for low-dimensional spaces, BO has recently been modified 
and used also for very large-dimensional spaces (up to about one thousand dimensions).

In this paper we consider a much simpler algorithm, called ``Reactive Affine Shaker'' (RAS)~\cite{BB2008}. The next sample is always generated with a uniform probability distribution
inside a parallelepiped (the ``box''). At each iteration, the form of the box is adapted during the search through an affine transformation,
based only on the point $\bsx$ position and on the success or failure in improving the function. The function values are therefore not used directly to modify the search area and to generate the next sample. The entire dimensionality is kept (no active subspaces like in \cite{papenmeier2022increasing}).

Despite its extreme simplicity and its use of only stochastic local search,
surprisingly the produced results are comparable to and not too far from the state-of-the-art results
of high-dimensional versions of BO, although with some more function evaluations.

An ablation study and an analysis of probability distribution of directions (improving steps and prevailing box orientation) in very large-dimensional spaces 
are conducted to understand more about the behavior of RAS and to assess the relative importance of the algorithmic building blocks
for the final results.

\end{abstract}

\section{Introduction}

All linear optimization problems are alike; each difficult nonlinear problem is difficult in its own way,
more so in high dimensions.
Each problem has different structural characteristics, and every specific instance of a problem
may present widely different ``fitness landscapes''. Furthermore, the local landscape of a single instance 
can depend critically on the position of the current point $\bsx$ in the input space and can therefore vary a lot during the search.
The so-called ``curse of dimensionality'' is well known. The complexity of finding an optimum -
especially with simple brute-force methods - tends to grow exponentially with the number of dimensions,
apart from rare special problems (like Linear Programming) or special instances of a problem.

Despite the negative theoretical worst-case results,
a concrete hope of finding improving solutions of practical interest derives from the fact that
many relevant real-world problems possess a high level of internal structure that can be exploited (e.g., the Big Valley hypothesis \cite{hains2011revisiting} or variations thereof)
and from the fact that the appropriate match between problem (or instance) characteristics
and algorithm configuration and parameters can be learned automatically via 
Machine Learning (ML). ML can be applied in an offline manner to adapt the meta-parameters to a problem
like in algorithm configuration
or online, by reacting to events occurring during a single search process on a specific instance 
(e.g., as advocated in Reactive Search Optimization - RSO~\cite{BBM2008thebook})

The growing availability of massive amounts of memory, 
starting from the eighties, opened new windows of opportunity for 
memory-based Intelligent Optimization techniques. The underlying assumption of a rich internal structure 
of most relevant optimization tasks makes techniques capable of 
\textit{gradually learning} that structure potentially more powerful 
and effective than memory-less techniques. 
 
This paper builds upon a previously proposed Reactive Search Optimization algorithm for global optimization
of multivariate functions of continuous variables called RAS~\cite{BB2006} (Reactive Affine Shaker). 
RAS is an adaptive search algorithm based only on point-wise 
function evaluations in a stochastic local (perturbative) search. 

The novel contributions of this paper are:
\begin{itemize}
     \item A qualitative and quantitative study of the evolution of the search box of RAS during the search,
          in particular for large-dimensional spaces
     \item An \textit{ablation study} of RAS to assess the relative contribution of its simple algorithmic building blocks 
     \item An experimental comparison of RAS on some very large dimensional problems previously
          solved with state-of-the-art Bayesian Optimization techniques
\end{itemize}

The remainder of this paper is organized as follows. In Section~\ref{sec:SoA} we discuss the state of the art on Bayesian optimization techniques; Section~\ref{sec:RAS} summarizes and motivates the proposed RAS heuristic. Experimental evaluations and comparisons are reported and discussed in Section~\ref{sec:experiments}. 
Appendix~\ref{sec:curse} discusses some issues that arise in high-dimensional optimization and provides additional motivation for some algorithmic choices in RAS.

\section{Bayesian Optimization for high-dimensional problems}
\label{sec:SoA}

In the context of the optimization of functions of continuous variables, we assume that the 
dominant computational cost is the evaluation of the 
function $f$ at sample points. This holds in many 
practical applications, e.g., when the evaluation of $f$ 
requires running a lengthy simulation, or even running an industrial 
plant and measuring the output. 

Bayesian optimization~\cite{frazier2018tutorial} is a sequential design strategy for global optimization of black-box functions
that are expensive to evaluate.
It builds a surrogate model for the objective from the sampled points, quantifies the uncertainty by using a Bayesian machine learning technique,
(Gaussian process regression), and then uses an \textit{acquisition function} defined from this surrogate to decide where to sample
the next point(s).
The acquisition function trades off exploration and exploitation to reduce the number of expensive function evaluations. 

Although based on Gaussian processes (without guarantees that concrete functions are appropriately modeled by GP) and expensive, BO is considered one of the most promising algorithms
for optimization of functions of limited dimensionality (usually not more than 10-20).

A series of breakthroughs have recently pushed the envelope of high-dimensional Bayesian optimization
for a wider adoption in science and engineering.
For the limited scope of this paper, we concentrate on two recent contributions~\cite{papenmeier2022increasing,eriksson2019scalable} 
and we refer to the contained bibliography for a review of the state of the art.

The research in~\cite{eriksson2019scalable} argues that the implicit
homogeneity of the global probabilistic models in BO tends to overemphasize exploration
in high-dimensional problems. 
The fact that search spaces grow faster than sampling budgets implies the presence of regions with large posterior uncertainty. For common
myopic acquisition functions, this results in an overemphasized exploration and a failure to exploit
promising areas.
To remedy, they
propose the TuRBO (Trust-region BO)  that considers a population of parallel and independent local models and performs a global allocation of samples across these models via an implicit multi-armed bandit
approach (the local models that are more promising get progressively more samples to evaluate).
The trust-region (TR) idea~\cite{yuan2000review} is that each local model can be trusted only in a region, typically a ball with a given radius around the current solution, a radius that is adapted during optimization.
In TuRBO the TR is a hyperrectangle centered at the best solution found so far $\bsx^*$.
The initial base side length $L$ is a portion of the range along each coordinate, the total volume of the TR is kept fixed, while the 
length along each  coordinate is rescaled according to its corresponding lengthscale $\lambda_i$ in the GP model.
The base side length $L$ is then adapted during the run, by keeping it sufficiently large so that the TR contains good solutions but small
enough to ensure that the local model is accurate within the TR. In detail, the TR is expanded after many consecutive ``successes''
and shrunk after many consecutive ``failures'', in a similar spirit as in \cite{nelder1965simplex}.

In the assumption of the existence of an \textit{active subspace}, with a projection matrix $T$ so that the value of the function
to optimize depends only on the projected $\bsx$ value ($F(\bsx) = g(T \bsx)$)
REMBO (Random embedding BO) \cite{wang2016bayesian} and HeSBO (Hashing-enhanced subspace BO) \cite{nayebi2019framework} try
to capture this active subspace by a randomly chosen linear subspace.
SaasBO~\cite{eriksson2021high} uses sparse priors on the GP length scales that is particularly effective if the active subspace is axis-aligned. Alebo \cite{letham2020re} uses a Mahalanobis kernel and
linear constraints on the acquisition function.

The work in \cite{papenmeier2022increasing} starts from the observation that methods for high-dimensional Bayesian optimization (HDBO) suffer from
degrading performance in spaces of growing dimension or risk failure
if some assumptions are not met. They propose a new algorithm (BAxUS -Bayesian optimization with adaptively expanding subspaces) that
uses the idea of nested low-dimensional random subspaces of growing dimensionality to adapt the space it optimizes
over to the problem.

\section{The RAS heuristic}
\label{sec:RAS}

We summarize and motivate in this section the RAS local search algorithm which is the method considered in this work.

\begin{figure}
%
%
%

\begin{tabular}{ll}
\hline
$f$ & Function to minimize\\
$\bsx$ & Initial point\\
$\calR$ & Search region\\
$\bsD$ & Current displacement\\
\hline
\end{tabular}

\begin{pseudo}
{\bf function} {\tt RAS} ($f$,  $\bsx$)\xl
\xa  $\calR$  $\leftarrow$ small isotropic set around  $\bsx$					\label{init R}\xl
\xb {\bf while} (local termination condition is not met)				\label{start repeat}\xl
\xb \xa Pick $\bsD\in\IR^d$ such that $\bsx+\bsD,\bsx-\bsD\in\calR$ \label{new delta}\xl
\xb \xb {\bf if} $f(\bsx+\bsD) < f(\bsx)$					\label{first shot}\xl
\xb \xb \xa  $\bsx$  $\leftarrow$  $\bsx$ +  $\bsD$;\xl
\xb \xb \xb Extend  $\calR$ along  $\bsD$\xl
\xb \xb \xc Center  $\calR$ on  $\bsx$\xl
\xb \xb {\bf else} {\bf if} $f(\bsx-\bsD) < f(\bsx)$				\label{second shot}\xl
\xb \xb \xa  $\bsx$  $\leftarrow$  $\bsx$ -  $\bsD$;\xl
\xb \xb \xb Extend  $\calR$ along  $\bsD$\xl
\xb \xb \xc Center  $\calR$ on  $\bsx$								\label{end second shot}\xl
\xb \xb {\bf else}\xl
\xb \xc \xn Reduce $\calR$ along  $\bsD$				\label{no shot}\xl
\xc {\bf return}  $\bsx$;
\end{pseudo}

\caption{The RAS algorithm} \label{fig:ras}
\end{figure}

The Reactive Affine Shaker Heuristic~\cite{BB2006}, RAS for short, 
is a self-tuning local search algorithm based on~\cite{solis-wets}. No prior knowledge is required 
on the function $f$ and only evaluations at arbitrary 
values of the independent variables are allowed. The RAS heuristic 
tries to rapidly move towards better objective values by maintaining 
and updating a ``search region'' $\calR$ around the current 
point $\bsx$.

The use of memory in RAS is limited: the entire previous history of 
the search (the trajectory of the generated sample points and the 
outcome of the evaluations) is summarized through the \textit{dynamic 
search region}, intended to zoom in on the promising areas where 
to find points better than the current record point.

RAS adapts the search 
region $\calR$ depending on the occurrence or lack of success during 
the last step. If a step in a certain direction is improving, then $\calR$ is expanded along that direction; it 
is reduced otherwise. Once a promising direction is 
found, the probability that subsequent steps will follow the same 
direction is increased, and the search will proceed more and more 
aggressively in that direction until bad results reduce its 
prevalence. The algorithm is outlined in Fig.~\ref{fig:ras}.

RAS starts with an isotropic search region centered around 
the initial point (line~\ref{init R}). Next, new sample points are 
 generated (line~\ref{new delta}). If the sample point 
$\bsx+\bsD$ yields a lower objective value (line~\ref{first shot} 
and following), then the current position is updated and $\calR$ is 
expanded along the direction of $\bsD$. To increase the probability 
of finding a better point, if $\bsx+\bsD$ does not lead to an 
improvement, also $\bsx-\bsD$ is tried (line~\ref{second shot} and 
following). If both points fail at improving $f$, then 
the search region is reduced along the direction of $\bsD$ 
(line~\ref{no shot}) and the current position is kept. The above steps are repeated until a local termination condition is 
verified. Common termination criteria are the number of 
iterations, the size of the search region, or a large 
number of iterations without further improvement.

\subsection{Implementation of the anisotropic search box}

The main rationale for the introduction of a non-isotropic search box lies in the increasing difficulty of finding improving directions, even for very smooth and regular functions, when displacements are chosen along uniformly random directions. These difficulties arise from the fact that in many dimensions most random directions tend to be orthogonal to the desired one (e.g., the function's gradient) and lead to worse values of the objective function. See Appendix~\ref{sec:curse} for an experimental discussion.

\begin{figure}[t]
    \centering
    \includegraphics[width=\hsize]{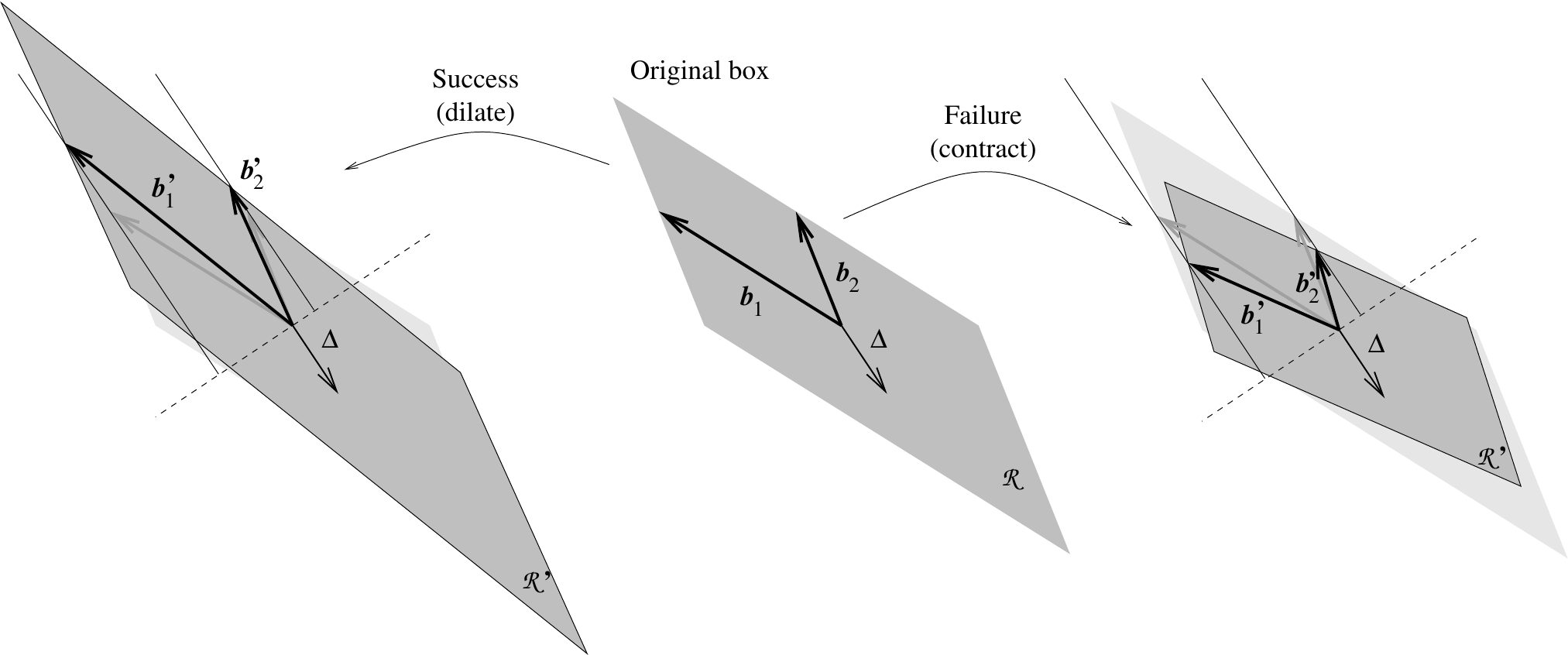}
    \caption{Evolution of the search box $\calR$ under transformation~\eqref{eq:affine}: from the current shape (center), the base vectors $\bsbe_i$ are dilated in the direction of $\bsD$ if the step succeeds (left), or contracted (right) if the step fails.}
    \label{fig:box_step}
\end{figure}

The search region $\calR$ is 
implemented as a box defined by $d$ independent base vectors 
($\bsbe_1\ldots \bsbe_d$), where $d$ is the number of dimensions of 
the search domain. The search region $\bsD$ is determined as a random linear combination of the base vector with independent coefficients uniformly distributed in $[-1,1]$:
$$
    \bsD=\sum_{i=1}^d r_i\bsbe_i, \qquad r_i\in[-1,1].
$$
Shape modifications are implemented as affine transformations of these vectors by contracting (if the step fails) or dilating (if the step succeeds) their components along the direction of $\bsD$ as shown in Fig.~\ref{fig:box_step}. Analytically, we want to add to every base vector $\bsbe_j$ a contribution in the direction of $\bsD$ proportional to the base vector's size:
$$
    (\rho-1)\frac{\bsD\cdot\bsbe_j}{\|\bsD\|^2}\bsD,
$$
where $\rho$ controls the amount of contraction ($0<\rho<1$, making the contribution opposite wrt to the projection of $\bsbe_j$ along $\bsD$) or dilation ($\rho>1$, making the contribution positive along the projection). This can be rewritten as the following linear transformation:
\begin{equation}\label{eq:affine}
    \forall j \quad \bsbe_j \leftarrow A\bsbe_j,\qquad\text{where}\qquad A=\bsi + (\rho-1)\frac{\bsD\bsD^T}{\|\bsD\|^2}.
\end{equation}

The RAS heuristic depends therefore on three parameters,
$$
0<\rho_\text{con}<1, \qquad \rho_\text{dil}>1, \qquad 0<\eta<1
$$
with the following meaning, assuming a search hyperinterval $[\text{min}_i,\text{max}_i]$, $i=1,\dots,d$:
\begin{itemize}
    \item the initial search box vectors $\bsbe_i$, $i=1,\dots,d$, are aligned along the domain axes, with length $\|\bsbe_i\|=\eta\cdot(\text{max}_i-\text{min}_i)$;
    \item the affine transformation factor in~\eqref{eq:affine} is $\rho=\rho_\text{dil}$ for dilations (upon improving step) and $\rho=\rho_\text{con}$ for contractions (upon double-shot failure).
\end{itemize}

\section{Experimental Results}
\label{sec:experiments}

The following parameter values are used in all experiments unless otherwise noted:
$$
    \eta=\frac15=0.2, \qquad \rho_\text{dil}=5, \qquad \rho_\text{con}=\frac1{\rho_\text{dil}}=0.2.
$$

Let's note that the dilation and contraction parameters are far from 1.
In a high-dimensionality setting, the affected direction $\bsD$ is, with high probability, almost perpendicular to all base vectors $\bsbe_i$, therefore the effect of milder parameters (i.e., close to 1) would be very limited.

\subsection{Behavior of the search box}

\begin{figure}[t]
    \centering
    \includegraphics[width=.49\hsize]{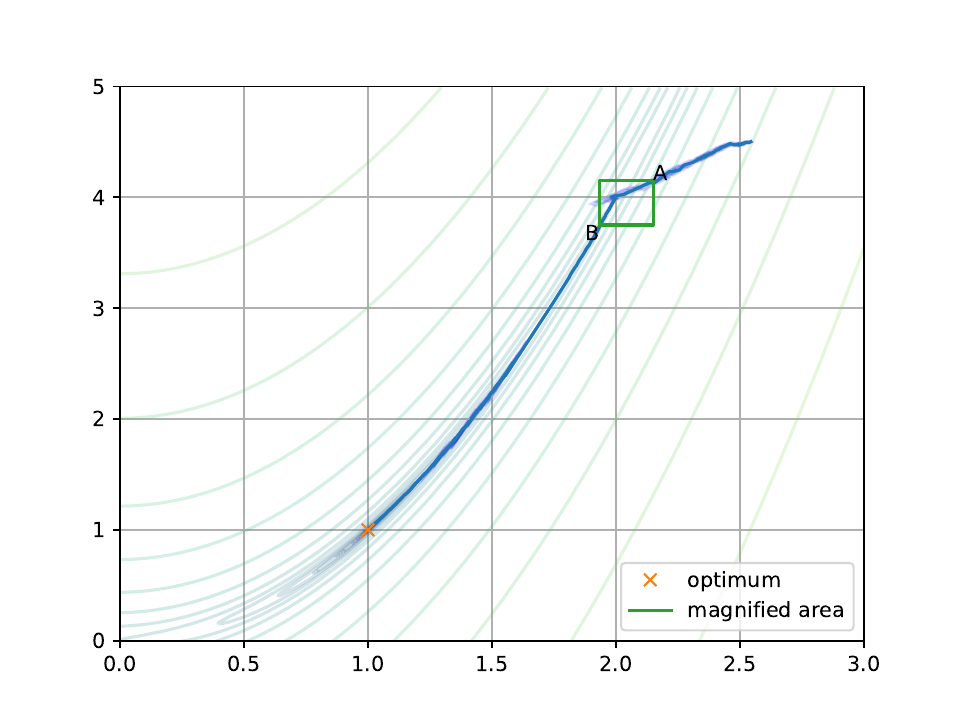}
    \includegraphics[width=.49\hsize]{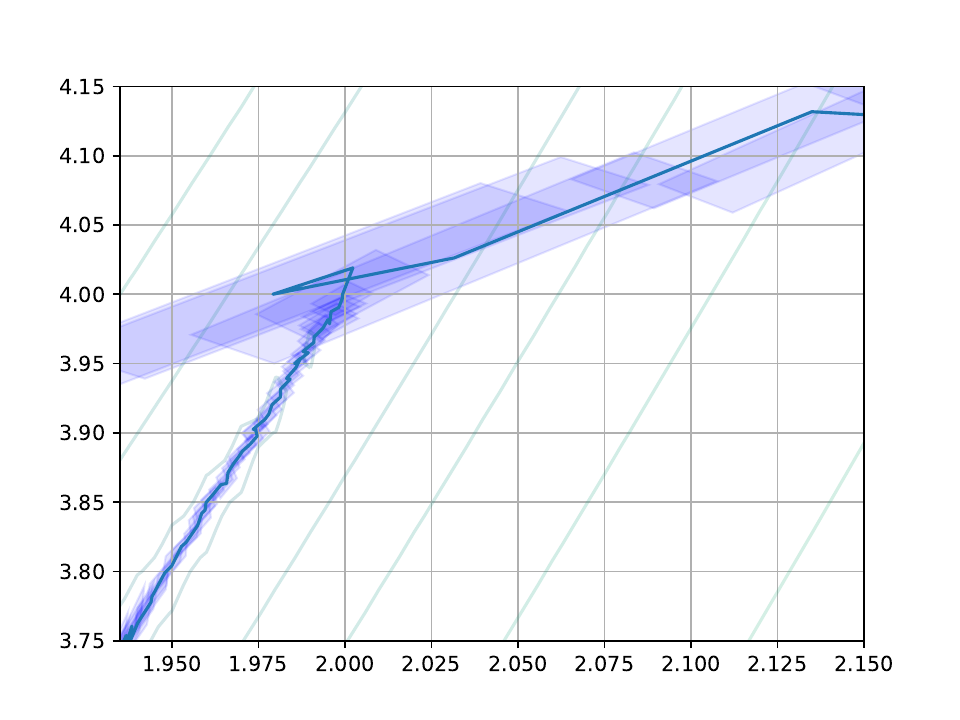}
    \caption{Evolution of the search box under transformation~\eqref{eq:affine}, 2D Rosenbrock function; on the right, a detailed view of the AB rectangle with the direction change.}
    \label{fig:box1}
\end{figure}

In Fig.~\ref{fig:box1} we show a sample RAS run on the 2D Rosenbrock function~\cite{rosenbrock1960automatic}, an unimodal function with a narrow, curved valley and with a selectable dimensionality. The global minimum is in $(1,1)$, and the search proceeds from right to left. The right-hand side plot shows the trajectory superimposed on the evolving search boxes. In the initial phase, the search proceeds downhill (we see the path crossing many function isolines); the search box is wide and slightly elongated in the search direction. Once the valley is reached around point $\bsx=(2.00,4.00)$, the trajectory proceeds with a narrower search box aligned with the valley direction.

\begin{figure}[t]
    \centering
    \begin{tabular}{cc}
        Rosenbrock (2D) & Paraboloid (100D)\\
        \includegraphics[width=.49\hsize]{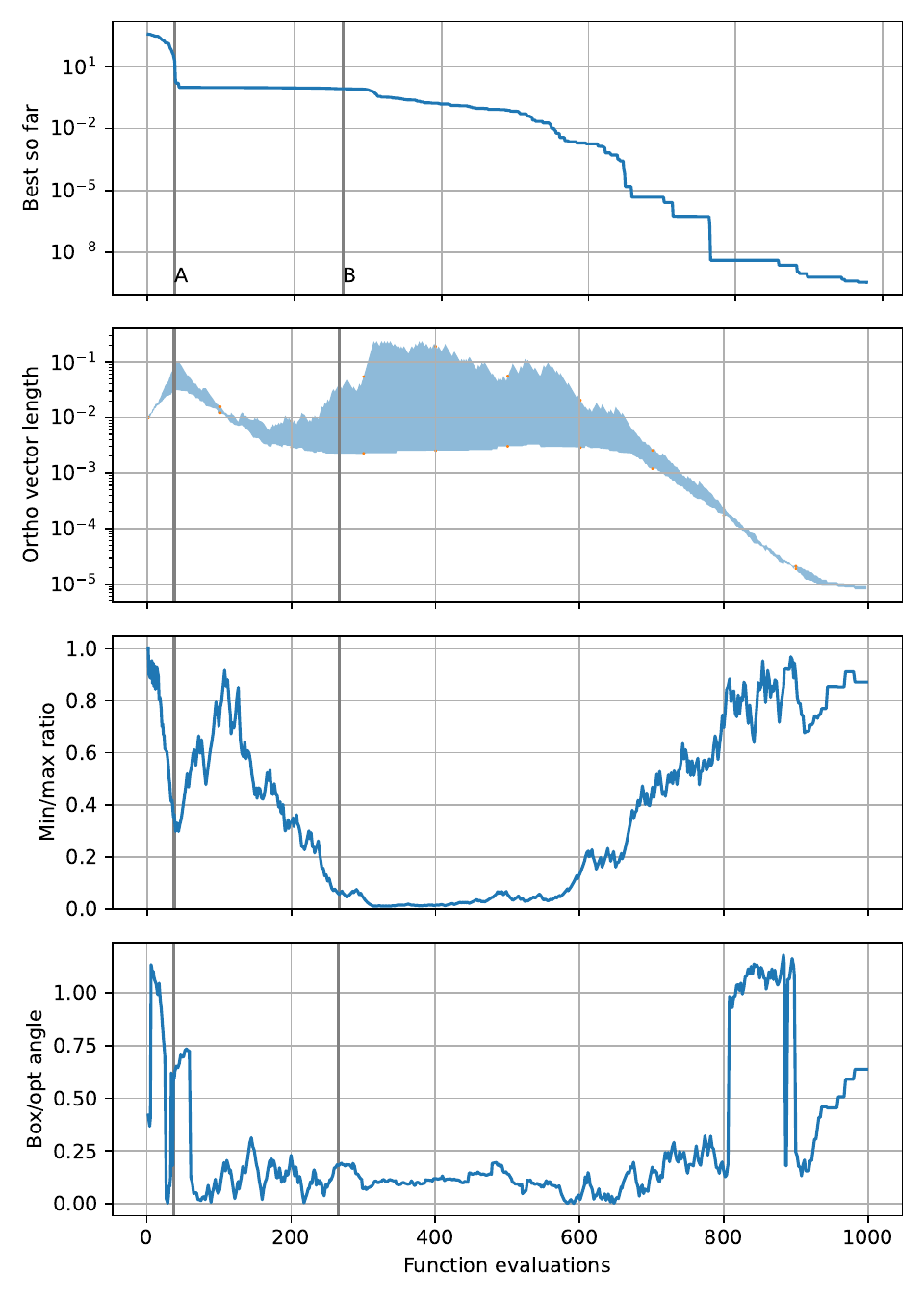} &
        \includegraphics[width=.49\hsize]{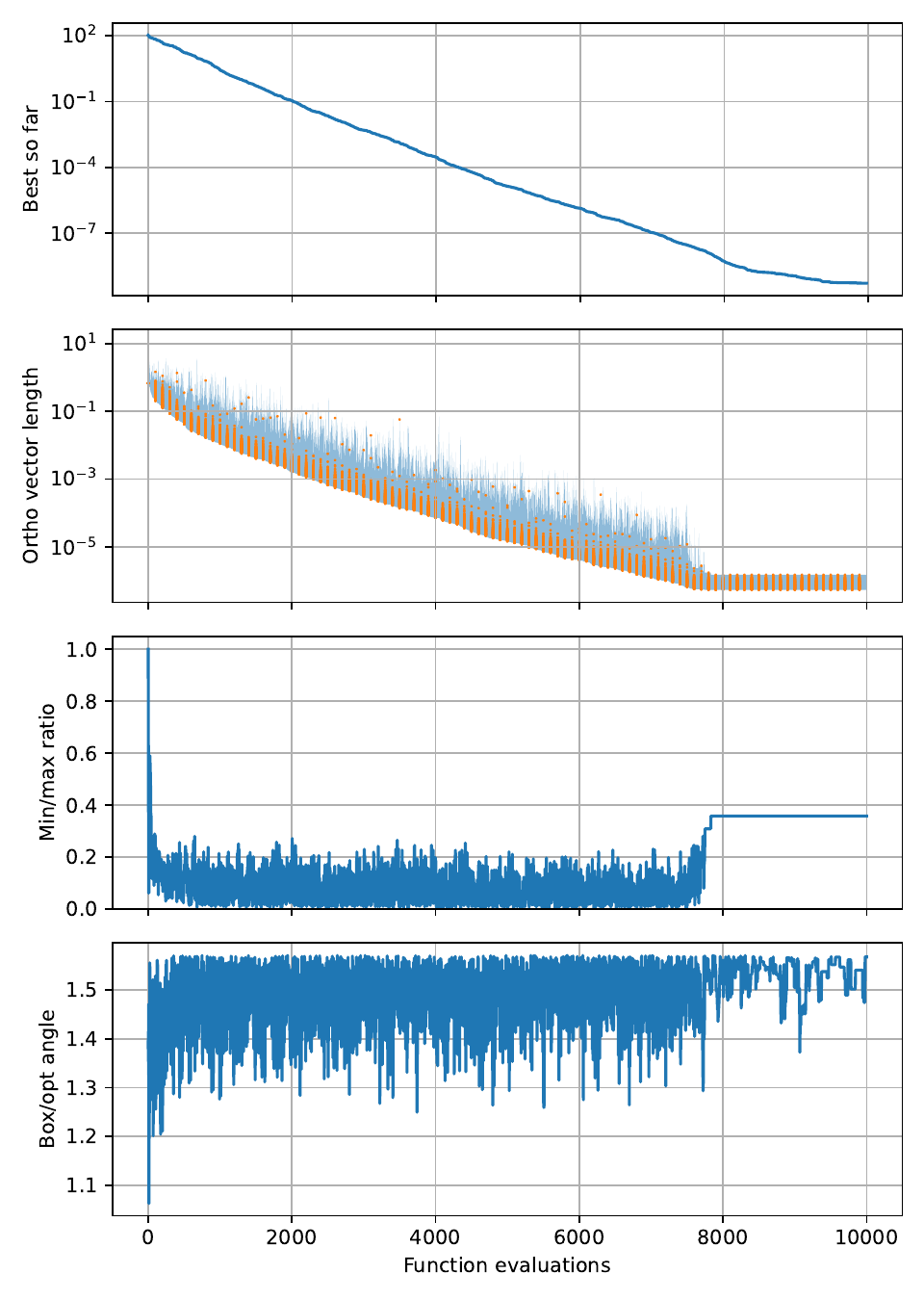}
    \end{tabular}
    \caption{Evolution of the search box under affine transformation wrt number of function evaluations for two test benchmarks, see the text for details.}
    \label{fig:box2}
\end{figure}

The analysis of the 2D Rosenbrock function search is expanded in the left column of Fig.~\ref{fig:box2}, where the evolution of the search is shown vs. the number of function evaluations. The top chart represents the best value found so far; the band delimited by the A and B markers corresponds to the portion of trajectory visible in the AB rectangle. The second chart from the top shows the size of the box, represented as a colored band spanning the size in the different directions; a wider band represents a more elongated search box. In the third chart, the evolution of the ratio between the minimum and maximum search box sizes is shown. Finally, the bottom chart displays the angle between the longest vector $\bsbe_i$ in the search box (i.e., the dominating search direction) and the direction from the current search position to the global optimum. The angle is given in radians, from $0$ (search box perfectly aiming at the optimal point) to $\pi/2\approx1.57$ (box aimed at orthogonal direction).

For the 2D Rosenbrock search, in the interval AB, we can see an initial fast improvement of the best value (the downhill phase). Once the valley is reached (very close to point A), the search box starts evolving in two ways: it becomes smaller (second row), because the shape of the valley causes many steps to fail, and more ``square'' (the ratio in the third plot approaches $1.0$). At some point, close to 100 evaluations, the search box starts elongating (ratio in third plot decreasing), and it points in the correct direction (fourth plot: oriented towards the global minimum). The dominant vector in the search box grows, encouraging steps in the right direction.

The right column of Fig.~\ref{fig:box2} shows the evolution of the search on a 100-dimensional paraboloid $f(\bsx)=\|\bsx\|^2$ on domain $\bsx\in[-1.5,1.5]^{100}$. Given the larger number of dimensions, more evaluations are required to achieve similar results. As the optimum point is approached, failures due to overshooting cause the box to shrink consistently (second row), while because of the regular shape, the box tends to stay quite elongated in a good direction (the third plot shows consistently small ratios between the shortest and the largest box vector). Although the bottom plot shows the dominant direction not to be aligned with the local optimum (and hence the gradient), this is  expected due to the high dimensionality of the search space and doesn't prevent the generation of improving moves.

\subsection{Comparison with state-of-the-art heuristics}

We compare the performance of RAS with a selection of heuristics covered in Section~\ref{sec:SoA} (TuRBO, SaasBO, Alebo, and HeSBO), the popular CMA-ES~\cite{hansen1996adapting} and random search~\cite{bergstra2012random} on the following six high-dimensionality benchmark problems with continuous variables on hyper-rectangular domains.\\
\textbf{Mopta08}~\cite{eriksson2021high} --- a 124-parameter vehicle design problem where the objective is the minimization of the vehicle's mass; the original problem is subject to several constraints transformed into soft penalties.\\
\textbf{SVM}~\cite{eriksson2021high} --- a Support Vector Machine training problem with 388 parameters.\\
\textbf{Branin2}~\cite{wang2016bayesian} --- the classical 2D test function with 3 global minimizers, embedded in a 500D space.\\
\textbf{Lasso-Hard}, \textbf{Lasso-High}~\cite{vsehic2022lassobench} --- two synthetic hard problems with~1000 and~300 parameters respectively from the \textsc{LassoBench} benchmark suite.\\
\textbf{Hartmann6}~\cite{wang2016bayesian} --- the classical 6D multimodal test function, embedded in a 500D space and rotated.

All functions are implemented in the publicly available BAxUS test suite\footnote{\url{https://baxus.papenmeier.io/}}~\cite{papenmeier2022increasing}, and all tests have been performed on the noiseless versions.

\begin{figure}[p]
    \centering
    \includegraphics[width=\hsize, trim={0 0 0 195mm}, clip]{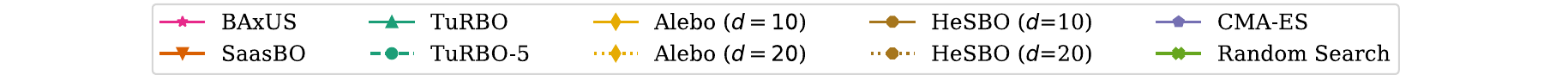}

    \vskip5mm
    \includegraphics[width=\hsize, trim={0 114.5mm 0 0}, clip]{lassomoptasvm_noiseless_svg-tex}
    \includegraphics[width=.32\hsize]{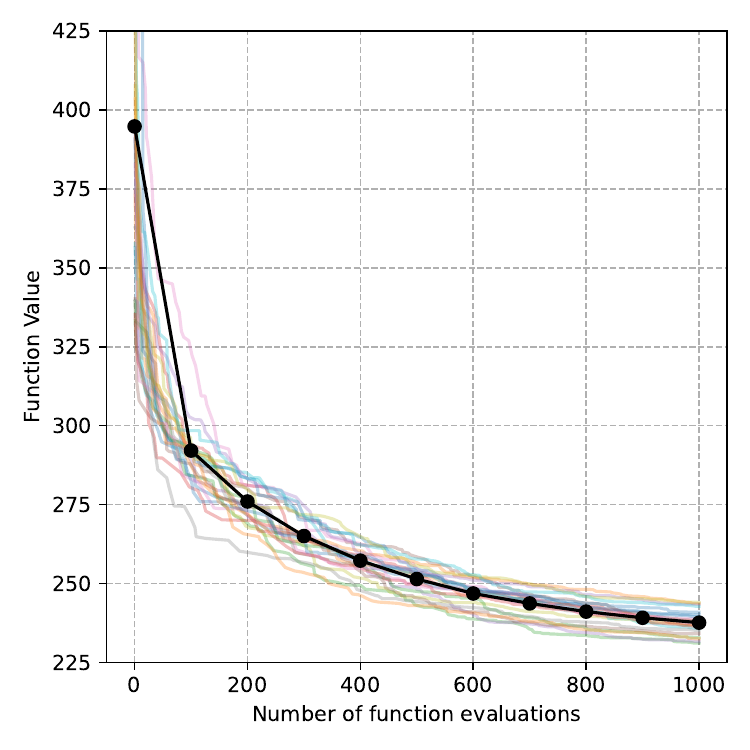}
    \includegraphics[width=.32\hsize]{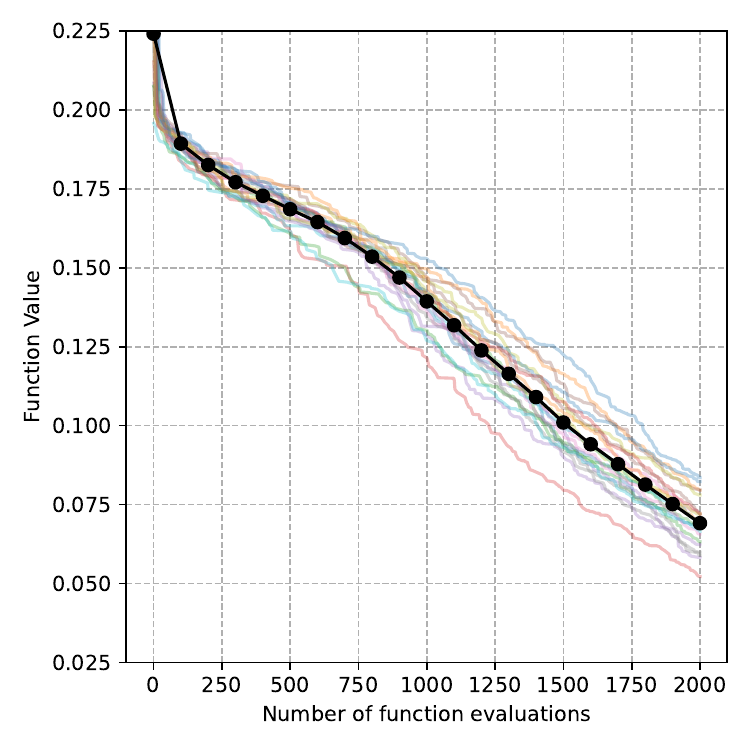}
    \includegraphics[width=.32\hsize]{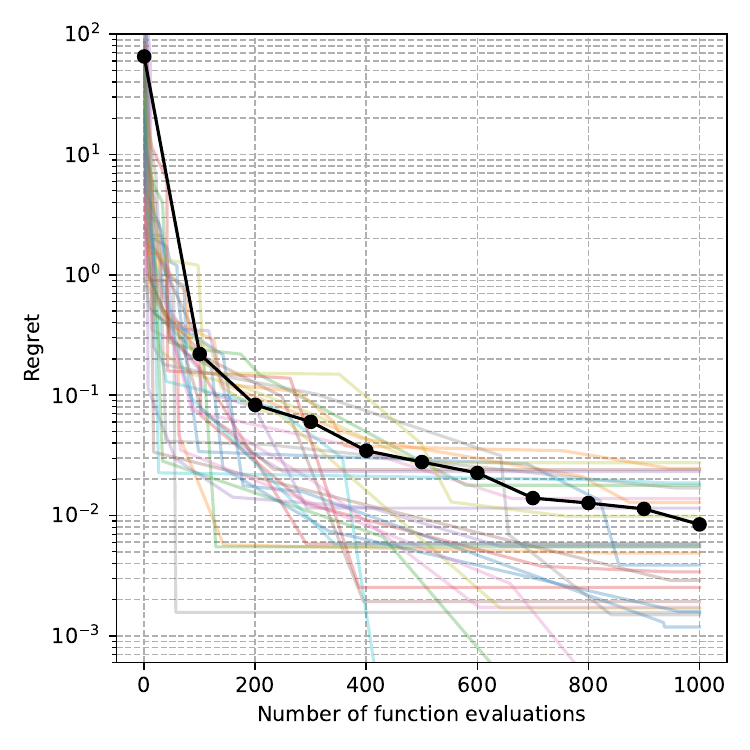}

    \vskip5mm
    \includegraphics[width=\hsize, trim={0 18mm 0 95mm}, clip]{lassomoptasvm_noiseless_svg-tex}
    \includegraphics[width=.32\hsize]{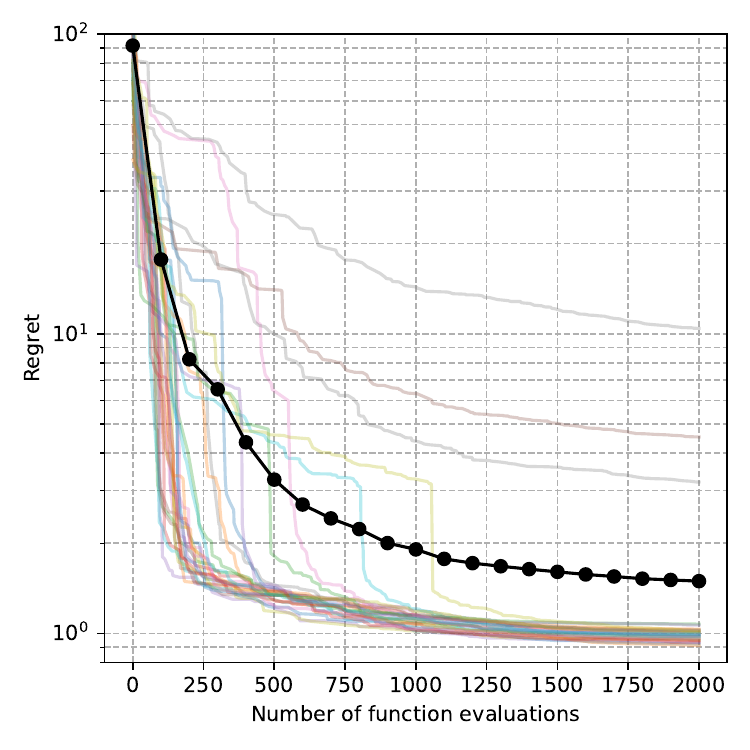}
    \includegraphics[width=.32\hsize]{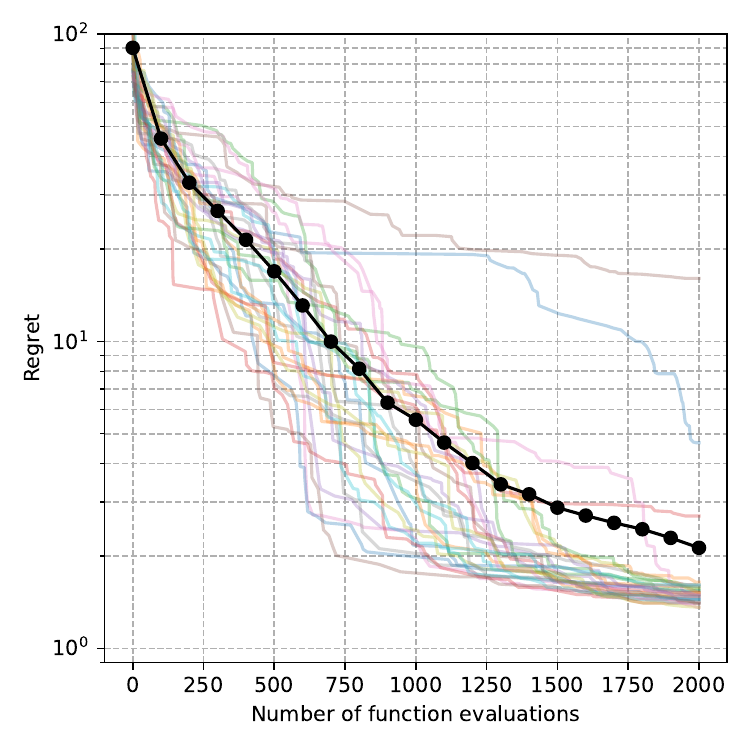}
    \includegraphics[width=.32\hsize]{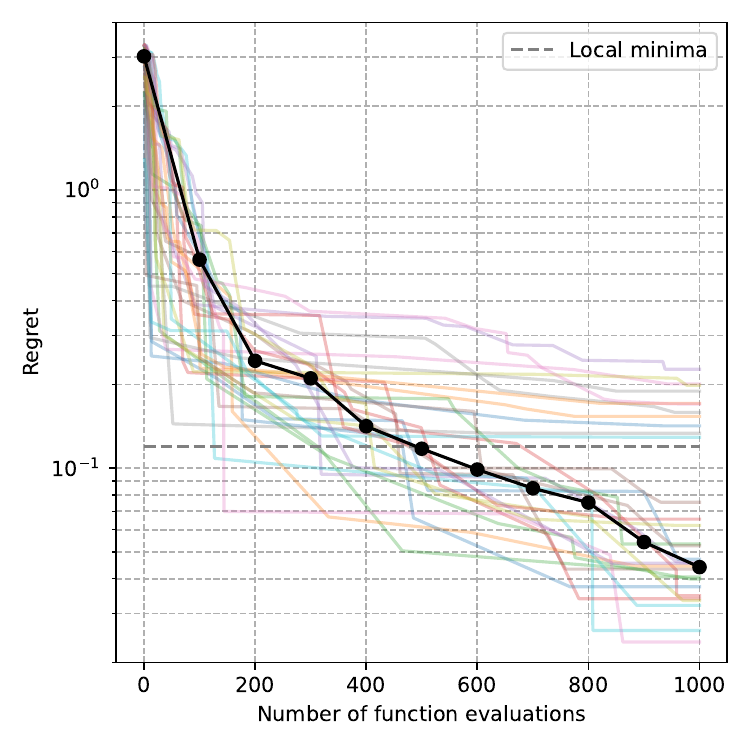}
    \caption{Test results on high-dimensionality functions. First and third row: TuRBO, SaasBO, Alebo, HeSBO, CMA-ES, random search (reprinted from~\cite{papenmeier2022increasing} with permission). Second and fourth row: results on the same functions from 30 runs of RAS, plotted on comparable scales.}
    \label{fig:plots}
\end{figure}

The results are shown in Fig.~\ref{fig:plots}. The top and third rows contain results from~\cite{papenmeier2022increasing}, while the second and bottom rows show the RAS results on the corresponding benchmark (names on top of each column) on comparable horizontal and vertical scales. For every function, all 30 RAS runs are shown as lightly colored lines, while the average is superimposed as a black dotted line.

RAS outperforms the Alebo, HeSBO, and SaasBO heuristics in all benchmarks with the exclusion of Branin2, where the identification of the 2D active subspace aligned with the domain axes kicks in very soon, while RAS keeps elongating and reducing the search box in irrelevant directions, and with the further exclusion of Hartmann6 for SaasBO. CMA-ES is outperformed in all cases except the Lasso-High benchmark.

On the other hand, the more sophisticated BAxUS and TuRBO heuristics generally behave better than RSA, with a significant advantage in the SVM, Branin2, and Lasso benchmarks. In the SVM and Lasso-High case, we note that RAS is still significantly decreasing at the end of the allotted number of evaluations, and further examination shows that $1.5\dots2\times$ evaluations are needed to achieve comparable results; in the three bottom benchmarks, RSA is further penalized by a few runs stuck in a slowly improving path (Lasso benchmarks) or in the attraction basin of a local minimum (Hartmann6, the dashed line represent the local minima).

Since RAS is a local optimization heuristic, the average line in the Hartmann6 plot only considers the trajectories in the global minimum attraction basin, corresponding to those that, at the 1000 evaluations mark, are below the dashed line. However, should the excluded trajectories be considered, the result would still outperform all heuristics except BAxUS and SaasBO.

\subsection{Ablation study}

To motivate the main algorithmic choices in the design of RAS, we performed a series of ablation studies in which specific aspects of the algorithm were silenced in order to assess their impact on test results.

\begin{figure}[t]
    \centering
    \includegraphics[width=\hsize]{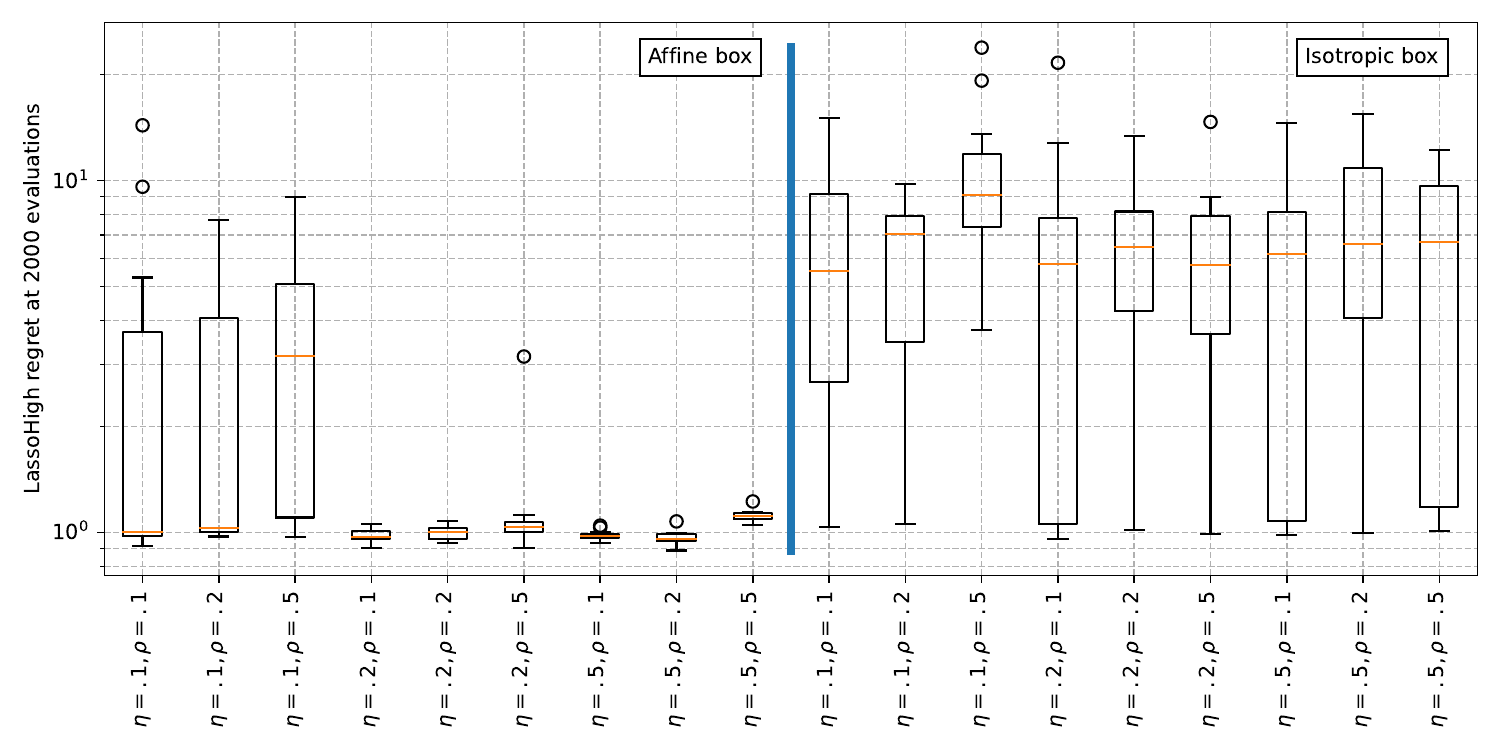}
    \caption{Effect of the search box affine transformation (left) compared with similarly-sized isotropic transformations (right).}
    \label{fig:isotropic box}
\end{figure}

\subsubsection{Isotropic search region}
To assess the impact of the affine transformation~\eqref{eq:affine} on search efficiency, we performed a series of tests in which the affine transformation~\eqref{eq:affine} has been replaced with a uniform (``isotropic'') resizing by the same factor: $\bsbe_i\leftarrow\rho\bsbe_i$. Fig.~\ref{fig:isotropic box} shows the results of a series of comparisons on the 300D Lasso-High function for different values of the search parameters $\eta$ (initial width of the search box wrt domain size) and contraction factor $\rho$ (the corresponding dilation factor is its reciprocal). Every boxplot collects the results of 30 tests. The results on the left-hand side of Fig.~\ref{fig:isotropic box} have been obtained with the unmodified RAS algorithm; the corresponding tests on the right-hand side of the same plot result from the isotropic version.\\
The effect of the affine transformation is very significant: with a high number of dimensions, the chance of finding an improving direction tends to reduce, and the anisotropic nature of transformation~\eqref{eq:affine} becomes the fundamental factor to guide the search without wasting function evaluations.

\begin{figure}[t]
    \centering
    \includegraphics[width=\hsize, trim={0 27mm 0 0}, clip]{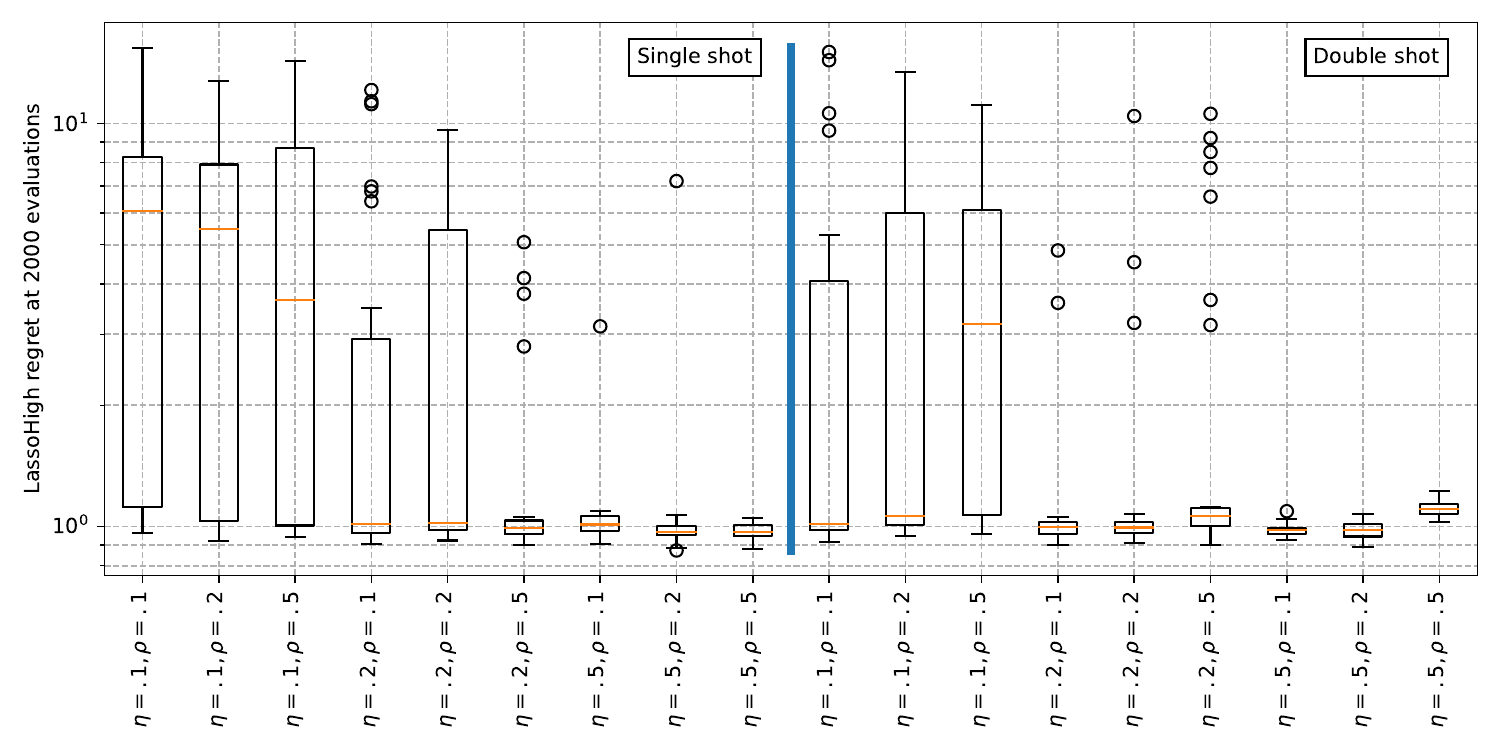}
    \includegraphics[width=\hsize]{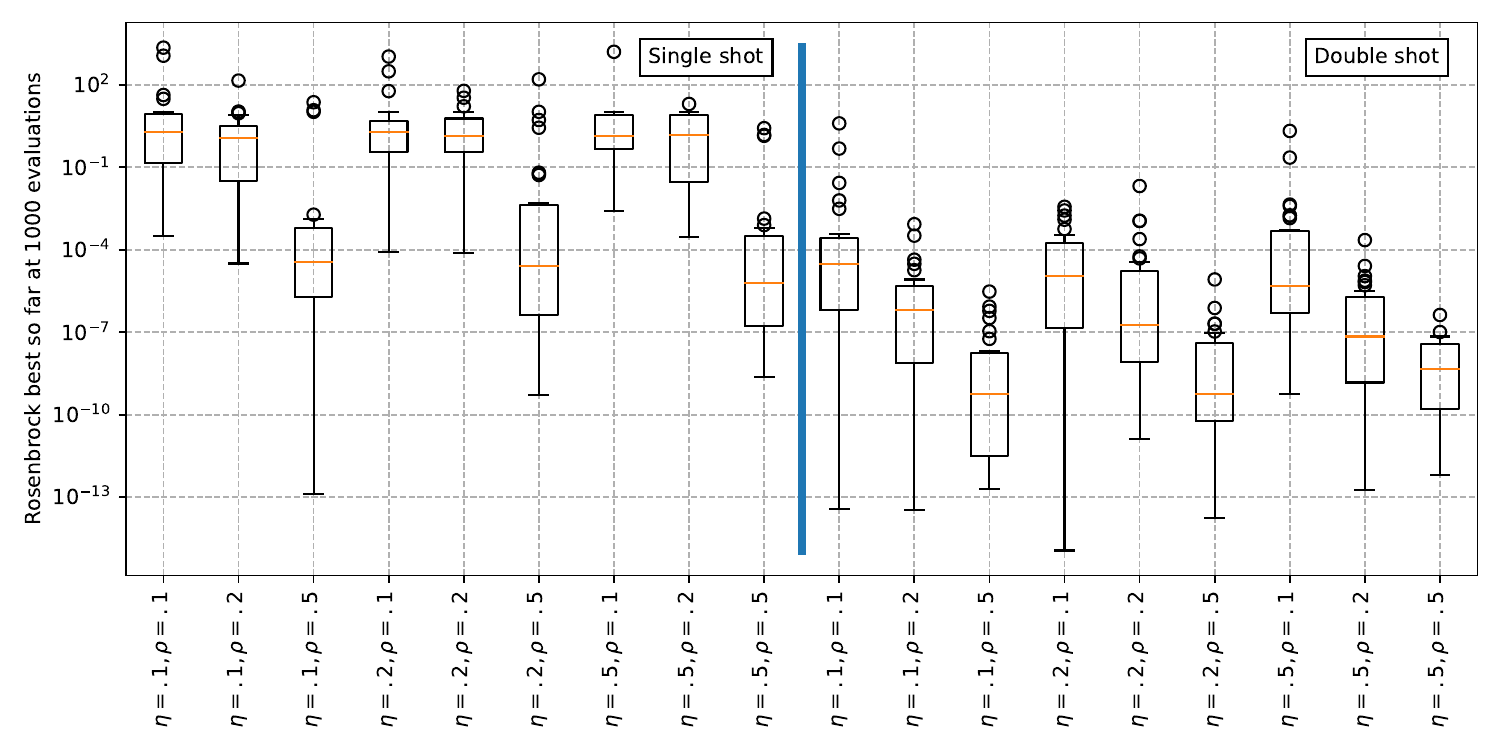}
    \caption{Effect of the double-shot strategy (top right) compared with a single-shot step (top left) on a high-dimensionality function. For reference, the bottom chart reports the same effect on the low-dimensionality Rosenbrock 2D function.}
    \label{fig:single shot}
\end{figure}

\subsubsection{Single-shot}
Another important component of the RAS search algorithm is the ``double-shot'' strategy: whenever a step in direction $\bsD$ fails, the opposite direction $-\bsD$ is tested. The rationale of this strategy is that a reasonably smooth function is locally approximable with a linear function, for which the double shot always works. However, it is expected that in high dimensionality functions, this advantage is reduced: because most randomly-generated search directions tend to be almost normal to the gradient, linear approximations only work in very small neighborhoods, and a second function evaluation before recalibrating the search box might be a waste of time.\\
In Fig.~\ref{fig:single shot}, we compare the full RAS algorithm (left-hand side of the plots) with a version where lines \ref{second shot}--\ref{end second shot} of Fig.~\ref{fig:ras} are removed (right-hand sides) for various parameter values and two different settings: the 300-dimensional Lasso-High function (top) and the 2D Rosenbrock function (bottom). While the advantage of the double-shot strategy is clear in the 2D experiment, its impact in many dimensions is quite interesting: although the results slightly worsen in some cases, its adoption seems to increase the algorithm's robustness concerning parameter changes; see in particular the results for $\eta=.2$, $\rho\in\{.1,.2\}$, where both the 3rd quartile and the inter-quartile range are greatly reduced in the double-shot case. These results call for further analysis on a wider parameter range and more diverse test functions.

\section{Conclusion}
\label{sec:conclusions}

The main conclusions of this study are:
\begin{itemize}
     \item RAS is surprisingly effective on the benchmark of very high dimensional problems.
          Despite its extreme simplicity, its use 
          of only success-failure information during the run (no consideration of the detailed function values $f(x_i)$ to build surrogate models), its stochastic local search nature (while BO is intended for global optimization)
     \item The results on the benchmark are comparable to those of the best BO algorithms
      (some more function evaluations are necessary but the difference is within a factor of two in most cases)
     \item RAS is superior w.r.t. many alternative techniques including popular Genetic Algorithms like  
     CMA-ES \cite{hansen1996adapting}
     although based on a very simple scheme (affine transformations of a search box and uniform probability for
     generating the next point)
     \item The ablation study confirms the extreme relevance of the affine transformation. 
     Instead of subspaces, a search box that is elongated and compressed with an affine transformation
     along successful or failure directions seems sufficient to reach interesting performance levels.
     \item Most of the problems in the considered benchmark do not seem to require global search methods (with the exception
     of Hartmann)
\end{itemize}

This preliminary investigation opens interesting questions about whether and when complex surrogate models
are critical for high-dimensional black-box function optimization. We plan to extend this research
through the integration of BO models on the box produced by RAS and to understand which (learnable) characteristics
of the functions are appropriate for the selection/configuration/tuning of different methods in the Intelligent Optimization paradigm.

At the same time, as already explored in low dimensions for the basic RAS scheme (like for the Repeated RAS of \cite{BBP2006}), 
we plan to investigate
parallel search streams to deal with multi-modal functions with more local optima, 
with a possible multi-armed bandit allocation strategy of with a Bayesian modeling of the potential of different starting points for initializing the local search.
A modification of RAS to deal with global optimization called $M-RAS$  via multiple parallel runs
and bandit-like allocation has been presented in 
\cite{BBP2006}.
$M-RAS$ is an extension of RAS in which promising starting
points for local search trails are suggested online by using Bayesian
Locally Weighted Regression.

The extension to combinatorial and mixed spaces (like the Bounce algorithm of \cite{papenmeier2023bounce}) is also on the stack.

\bibliographystyle{splncs04}
\bibliography{biblio}

\begin{thebibliography}{10}
\providecommand{\url}[1]{\texttt{#1}}
\providecommand{\urlprefix}{URL }
\providecommand{\doi}[1]{https://doi.org/#1}

\bibitem{BBM2008thebook}
Battiti, R., Brunato, M., Mascia, F.: Reactive Search and Intelligent
  Optimization, Operations research/Computer Science Interfaces, vol.~45.
  Springer Verlag (2008)

\bibitem{bergstra2012random}
Bergstra, J., Bengio, Y.: Random search for hyper-parameter optimization.
  Journal of machine learning research  \textbf{13}(2) (2012)

\bibitem{BB2008}
Brunato, M., Battiti, R.: {RASH}: A self-adaptive random search method. In:
  Cotta, C., Sevaux, M., S{\"o}rensen, K. (eds.) Adaptive and Multilevel
  Metaheuristics, Studies in Computational Intelligence, vol.~136. Springer
  (2008)

\bibitem{BB2006}
Brunato, M., Battiti, R.: The reactive affine shaker: a building block for
  minimizing functions of continuous variables. Tech. Rep. DIT-06-012,
  Universit{\`a} di Trento (Feb 2006)

\bibitem{BBP2006}
Brunato, M., Battiti, R., Pasupuleti, S.: A memory-based rash optimizer. In:
  Geffner, A.F.R.H.H. (ed.) Proceedings of AAAI-06 workshop on Heuristic
  Search, Memory Based Heuristics and Their applications. pp. 45--51. Boston,
  Mass. (2006), iSBN 978-1-57735-290-7

\bibitem{eriksson2021high}
Eriksson, D., Jankowiak, M.: High-dimensional bayesian optimization with sparse
  axis-aligned subspaces. In: Uncertainty in Artificial Intelligence. pp.
  493--503. PMLR (2021)

\bibitem{eriksson2019scalable}
Eriksson, D., Pearce, M., Gardner, J., Turner, R.D., Poloczek, M.: Scalable
  global optimization via local bayesian optimization. Advances in neural
  information processing systems  \textbf{32} (2019)

\bibitem{frazier2018tutorial}
Frazier, P.I.: A tutorial on bayesian optimization. arXiv preprint
  arXiv:1807.02811  (2018)

\bibitem{hains2011revisiting}
Hains, D.R., Whitley, L.D., Howe, A.E.: Revisiting the big valley search space
  structure in the tsp. Journal of the Operational Research Society
  \textbf{62}(2),  305--312 (2011)

\bibitem{hansen1996adapting}
Hansen, N., Ostermeier, A.: Adapting arbitrary normal mutation distributions in
  evolution strategies: The covariance matrix adaptation. In: Proceedings of
  IEEE international conference on evolutionary computation. pp. 312--317. IEEE
  (1996)

\bibitem{letham2020re}
Letham, B., Calandra, R., Rai, A., Bakshy, E.: Re-examining linear embeddings
  for high-dimensional bayesian optimization. Advances in neural information
  processing systems  \textbf{33},  1546--1558 (2020)

\bibitem{nayebi2019framework}
Nayebi, A., Munteanu, A., Poloczek, M.: A framework for bayesian optimization
  in embedded subspaces. In: International Conference on Machine Learning. pp.
  4752--4761. PMLR (2019)

\bibitem{nelder1965simplex}
Nelder, J.A., Mead, R.: A simplex method for function minimization. The
  computer journal  \textbf{7}(4),  308--313 (1965)

\bibitem{papenmeier2022increasing}
Papenmeier, L., Nardi, L., Poloczek, M.: Increasing the scope as you learn:
  Adaptive bayesian optimization in nested subspaces. Advances in Neural
  Information Processing Systems  \textbf{35},  11586--11601 (2022)

\bibitem{papenmeier2023bounce}
Papenmeier, L., Nardi, L., Poloczek, M.: Bounce: reliable high-dimensional
  bayesian optimization for combinatorial and mixed spaces. Advances in Neural
  Information Processing Systems  \textbf{36},  1764--1793 (2023)

\bibitem{rosenbrock1960automatic}
Rosenbrock, H.: An automatic method for finding the greatest or least value of
  a function. The computer journal  \textbf{3}(3),  175--184 (1960)

\bibitem{vsehic2022lassobench}
{\v{S}}ehi{\'c}, K., Gramfort, A., Salmon, J., Nardi, L.: Lassobench: A
  high-dimensional hyperparameter optimization benchmark suite for lasso. In:
  International Conference on Automated Machine Learning. pp.~2--1. PMLR (2022)

\bibitem{solis-wets}
Solis, F.J., Wets, R.J.B.: Minimization by random search techniques.
  Mathematics of Operations Research  \textbf{6}(1),  19--30 (1981)

\bibitem{wang2016bayesian}
Wang, Z., Hutter, F., Zoghi, M., Matheson, D., De~Feitas, N.: Bayesian
  optimization in a billion dimensions via random embeddings. Journal of
  Artificial Intelligence Research  \textbf{55},  361--387 (2016)

\bibitem{yuan2000review}
Yuan, Y.x.: A review of trust region algorithms for optimization. In: Iciam.
  vol.~99-1, pp. 271--282 (2000)

\end{thebibliography}

\clearpage

\appendix

\section{Supplemental material --- Finding improving directions in high dimensions}
\label{sec:curse}

The main rationale behind the proposal of RAS for high-dimensional functions lies in the increasing difficulty of finding improving directions even for very smooth and regular functions, if such directions are uniformly sampled. This difficulty rises from two main factors, discussed in the following subsections.

\subsection{Random directions tend to be mutually orthogonal}

\begin{table}[b]
    \caption{Average angle between random vectors in $d$ dimensions (top row) and expected number of double-shot successes for a random displacement depending on the search box radius $r(B')$ relative to the average curvature of the function's isosurface.}
    \label{tab:curvature}
    \centering
    \begin{tabular}{l|rrrrrrrrr}
        & \multicolumn{9}{c}{Dimension $d$}\\
        & \multicolumn{1}{c}{1}
        & \multicolumn{1}{c}{2}
        & \multicolumn{1}{c}{3}
        & \multicolumn{1}{c}{5}
        & \multicolumn{1}{c}{10}
        & \multicolumn{1}{c}{50}
        & \multicolumn{1}{c}{100}
        & \multicolumn{1}{c}{500}
        & \multicolumn{1}{c}{1000}\\
        \hline
        &&&&&&&&&\\[-1em]
        Avg angle $\bar\theta_d$ (degrees)	&0.00	&45.00	&57.30	&66.85	&74.64	&83.46	&85.40	&87.95	&88.55\\
        \hline
        $r_{B'}/r_B=1.000$	&1.00	&0.78	&0.62	&0.42	&0.16	&0.00	&0.00	&0.00	&0.00\\
        $r_{B'}/r_B=0.500$	&1.00	&0.89	&0.81	&0.69	&0.50	&0.08	&0.01	&0.00	&0.00\\
        $r_{B'}/r_B=0.100$	&1.00	&0.98	&0.96	&0.94	&0.89	&0.73	&0.62	&0.26	&0.11\\
        $r_{B'}/r_B=0.050$	&1.00	&0.99	&0.98	&0.97	&0.95	&0.86	&0.81	&0.58	&0.43\\
        $r_{B'}/r_B=0.010$	&1.00	&1.00	&1.00	&0.99	&0.99	&0.97	&0.96	&0.91	&0.87\\
        $r_{B'}/r_B=0.005$	&1.00	&1.00	&1.00	&1.00	&0.99	&0.99	&0.98	&0.95	&0.94\\
        $r_{B'}/r_B=0.001$	&1.00	&1.00	&1.00	&1.00	&1.00	&1.00	&1.00	&0.99	&0.99\\
    \end{tabular}
\end{table}

Following~\cite{BB2008}, the average angle $\bar\theta_d$ in radians between two random directions in $d\ge2$ dimensions is
$$
    \bar\theta_d=\frac{J_{d-2}}{I_{d-2}},
$$
where the numerator and denominator are recursively defined as
$$
    I_d=\begin{cases}
        \frac\pi2&\text{if $d=0$}\\
        1&\text{if $d=1$}\\
        \frac{(d-1)I_{d-2}}d&\text{if $d>1$,}
    \end{cases}
    \qquad
    J_d=\begin{cases}
        \frac{\pi^2}8&\text{if $d=0$}\\
        1&\text{if $d=1$}\\
        \frac{(d-1)J_{d-2}}d+\frac1{d^2}&\text{if $d>1$.}
    \end{cases}
$$
The first row of Table~\ref{tab:curvature} reports the values of $\bar\theta$ computed for an increasing number of dimensions.

For $d=2$ the average angle is $45\degree$, but for higher and higher dimensionalities the average tends to $90\degree$. In the simplifying hypothesis that the function's gradient is unknown and that we generate a move in a random direction, this means that in high dimensions there is a progressively smaller probability of moving along the gradient.

\subsection{Improving moves tend to be scarce}

\begin{figure}[tbp]
    \centering
    \includegraphics[width=.6\hsize]{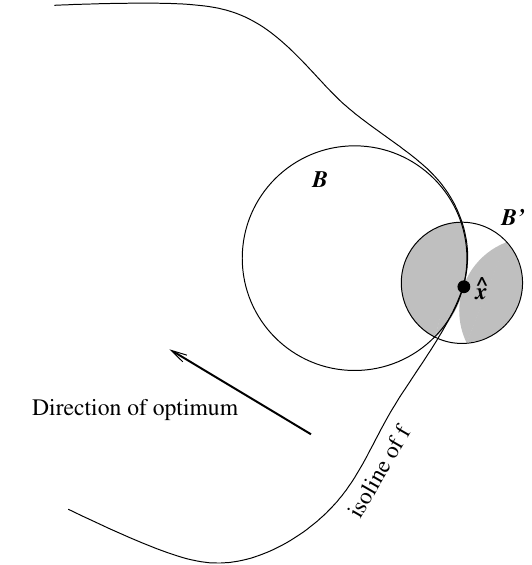}
    \caption{Estimation of the probability of success for the double-shot strategy with simplifying hypotheses (locally smooth function with isolines locally approximable with a spherical surface $B$, round search box $B'$). The grey area defines a successful displacement (considering the double shot strategy); the ratio between the grey area and the area of $B'$ gives the success probability of the double-shot strategy.}
    \label{fig:spheres}
\end{figure}

The effectivenes of the double-shot strategy in the limit for small search boxes was proved in~\cite{BB2008}. However, given a desired success probability, as the number of search-space dimensions grows the required search box size might become impractically small.
Consider the situation described in Fig.~\ref{fig:spheres}. Let $\hat\bsx\in\IR^d$ be the current point during a local search for the minimum of function $f:D\to\IR$ on domain $D\subseteq\IR^d$.

We define the isosurface of function $f$ at point $\hat\bsx$ to be the locus of all points $\bsx$ in its domain with the same function value:
$$
    \text{iso}_f(\hat\bsx) = f^{-1}(f(\hat\bsx)) = \{\bsx:f(\bsx)=f(\hat\bsx)\}.
$$
The isosurface has dimension $d-1$ and is embedded in $\IR^d$. It locally divides the domain into two sides, the ``good'' one where $f(\bsx)<f(\hat\bsx)$ (left side of the isoline in Fig.~\ref{fig:spheres}) and the ``bad'' one where $f(\bsx)>f(\hat\bsx)$. If $f$ is locally convex and smooth enough, we can approximate the curvature radius of the iso-surface at point $\hat\bsx$ by considering the radius of the largest sphere $B$ tangent to the iso-surface in $\hat\bsx$ and lying on the same side of it\footnote{For $d>2$, we actually have different curvature radii on different directions; let us consider the smallest one.}.

Suppose now that the local move consists of generating a random displacement $\bsD$ in a ball $B'$ centered on $\hat\bsx$. The move succeeds if $\hat\bsx+\bsD$ falls inside $B$ (as a local approximation of the ``interior'' of the isosurface) or, due to the double-shot strategy, in the symmetrically opposite direction. The probability of success is therefore given by the ratio between the area of the greyed-out portion of $B'$ in Fig.~\ref{fig:spheres} and the full area of $B'$, and it depends on the ratio between the radius of $B'$ and the radius of $B$. We can consider $r_B$ as a normalization factor for this discussion.

Table~\ref{tab:curvature} shows the experimental success probability of the double-shot strategy at different dimensions $d$ and for different ratios between the search area radius $r_{B'}$ and the isosurface's curvature radius $r_B$. In particular, we can see that for high dimensionalities the success probability of the double-shot strategy is negligible and only becomes significant for small search box radii, forcing a heuristic based on an isotropic search box to adopt small displacements that slow the progress down.

The RAS heuristic tries to compensate this problem by introducing an anisotropic search region that elongates in an orthogonal direction with respect to the isosurface.

\end{document}